\documentclass[12pt]{article}
\usepackage[dvipsnames]{xcolor}
\usepackage{amsmath}
\usepackage{setspace}
\usepackage{bm}
\usepackage{bbm}
\usepackage{amssymb}
\usepackage{indentfirst}
\usepackage[superscript]{cite}
\usepackage{geometry}
\usepackage{mathtools}
\usepackage{tocloft}

\newcommand*{\citena}[1]{%
\begingroup
[\color{Green}
\romannumeral-`\x 
\setcitestyle{numbers}%
\cite{#1}%
\endgroup
]\ignorespacesafterend
}

\newcommand*{\citesup}[1]{%
\begingroup
\color{Green}
\cite{#1}%
\endgroup
\ignorespacesafterend
}

\newcommand*{\eqrefe}[1]{%
\begingroup
(\color{BrickRed}
\romannumeral-`\x 
\setcitestyle{numbers}%
\ref{eq:#1}%
\endgroup
)\ignorespacesafterend
}

\newtagform{Tags}[\textcolor{BrickRed}]{\color{Black}(}{)}

\DeclarePairedDelimiter\abs{\lvert}{\rvert}%

\usepackage[super]{natbib}
\setcitestyle{super}

\newcommand{\ii}{\bm{i}}

\DeclareMathOperator{\csch}{csch}
\DeclarePairedDelimiter{\floor}{\lfloor}{\rfloor}

\begin{document}
\title{A Reformulation of the Riemann Hypothesis}
\date{November 28, 2020}
\author{Jose Risomar Sousa}
\maketitle
\usetagform{Tags}

\begin{abstract}
We present some novelties on the Riemann zeta function. Using an extended formula created for the polylogarithm in a previous paper, $\mathrm{Li}_{k}(e^{z})$, the zeta function's Dirichlet series is analytically continued from $\Re(k)>1$ to the right half-plane, $\Re(k)>0$, by means of the Dirichlet eta function. More strikingly, we offer a reformulation of the Riemann hypothesis through a zeta's cousin, $\varphi(k)$, a pole-free function defined on the entire complex plane whose non-trivial zeros coincide with those of the zeta function.
\end{abstract}


\section{Introduction}
The Riemann Hypothesis is a long-standing problem in math, which involves the zeros of the analytic continuation of its most famous Dirichlet series, the zeta function. This Dirichlet series, along with its analytic continuation, constitutes a so-called $L$-function, whose zeros encode information about the location of the prime numbers. Riemann provided insight into this connection through his unnecessarily convoluted prime counting functions\citesup{Prime}.\\

The zeta function as a Dirichlet series is given by,
\begin{equation} \nonumber
\zeta(k)=\sum_{j=1}^{\infty}\frac{1}{j^k} \text{,}
\end{equation}
and throughout here we use $k$ for the variable instead of the usual $s$, to keep the same notation used in previous papers released on generalized harmonic numbers and progressions and interconnected subjects.\\ 

This series only converges for $\Re(k)>1$, but it can be analytically continued to the whole complex plane except $k=1$, the pole of the function. For the purpose of analyzing the zeros of the zeta function, we produce its analytic continuation to the right complex half-plane only, $\Re{(k)>0}$, by means of the alternating zeta function, known as the Dirichlet eta function, $\eta{(k)}$. It is a well known fact that all the non-trivial zeros of $\zeta(k)$ lie on the critical strip ($0<\Re{(k)<1}$). The Riemann hypothesis is then the conjecture that all such zeros have $\Re{(k)}=1/2$.\\

The exposition starts from the extended formula we found for the polylogarithm function, discussed in paper \citena{AC}. The polylogarithm is a generalization of the zeta function, and it has the advantage of having the Dirichlet eta function as a particular case.\\

The idea is to come up with a simplified equation on a complex variable, which preserves the non-trivial zeros of the zeta function, and then turn it into a system of two equations on two real variables (while also simplifying the convoluted expressions obtained by means of simple transformations). At the end of process, the system of equations is turned back into a single equation on a complex variable.

\section{The polylogarithm, $\mathrm{Li}_{k}(e^{z})$}
As seen in \citena{AC}, the following expression for the polylogarithm holds for all complex $k$ with positive real part, $\Re{(k)}>0$, and all complex $z$ (except $z$ such that $\Re{(z)}>=0$ and $\abs{\Im{(z)}}>2\pi$ -- though for $\abs{\Im{(z)}}=2\pi$ one must have $\Re{(k)}>1$),
\begin{multline} \label{eq:Polylog_form}
\mathrm{Li}_{k}(e^{z})=-\frac{z^{k}}{2\,\Gamma(k+1)}-\frac{z^{k-1}\left(1+\log{(-z)}\right)}{\Gamma(k)}\\
-\frac{z^{k}}{2\,\Gamma(k)}\int_{0}^{1}(1-u)^{k-1}\coth{\frac{z u}{2}}+\frac{2\left(1-z^{-k}(z-\log{(1-u)})^k\right)}{k\,u^2}\,du
\end{multline}\\
\indent From this formula, two different formulae for $\zeta(k)$ can be derived, using $z=0$ or $z=2\pi\ii$, but both are only valid when $\Re{(k)}>1$. Using $z=2\pi\ii$ is effortless, we just need to replace $z$ with $2\pi\ii$ in the above. Using $z=0$ is not as direct, we need to take the limit of $H_k(n)$\footnote{A formula derived from the partial sums of $\mathrm{Li}_{k}(e^{z})$, as explained in \citena{AC}.} as $n$ tends to infinity,
\begin{equation} \nonumber
\zeta(k)=\frac{1}{\Gamma(k+1)}\int_{0}^{1}\frac{\left(-\log{u}\right)^k}{(1-u)^2}\,du \text{}
\end{equation}

\subsection{The analytic continuation of $\zeta{(k)}$}
The analytic continuation of the zeta function to the right complex half-plane can be achieved using the Dirichlet eta function, as below,
\begin{equation} \label{eq:eta_def}
\zeta(k)=\frac{1}{1-2^{1-k}}\eta(k)=\frac{1}{1-2^{1-k}}\sum_{j=1}^{\infty}\frac{(-1)^{j+1}}{j^{k}} \text{,}
\end{equation}
\noindent which is valid in the only region that matters to the zeta function's non-trivial zeros, the critical strip. The exception to this mapping are the zeros of $1-2^{1-k}$, which are also zeros of $\eta(k)$, thus yielding an undefined product.\\

For the purpose of studying the zeros of the zeta function, we can focus only on $\eta(k)$ and ignore its multiplier. The below should hold whenever $\Re{(k)}>0$,
\begin{multline} \label{eq:eta_formula}
\eta{(k)}=-\mathrm{Li}_{k}(-1)=\frac{(\ii\pi)^{k}}{2\,\Gamma(k+1)}+\frac{(\ii\pi)^{k-1}\left(1+\log{(-\ii\pi)}\right)}{\Gamma(k)}\\
+\frac{(\ii\pi)^{k}}{2\,\Gamma(k)}\int_{0}^{1}-\ii(1-u)^{k-1}\cot{\frac{\pi u}{2}}+\frac{2\left(1-\pi^{-k}(\pi+\ii\log{(1-u)})^k\right)}{k\,u^2}\,du
\end{multline}

\section{Simplifying the problem}
The first step is to choose a suitable transform of $\eta{(k)}$ to simplify the calculation, such as,
\begin{equation} \label{eq:h(k)_eta_rel}
h(k)=2\,\Gamma(k+1)(\ii\pi)^{-k}\,\eta(k) \text{,}
\end{equation}
\noindent since $\eta{(k)}=0$ if and only if $h(k)=0$. This gives the below equation,
\begin{multline} \nonumber
-1+\frac{2\ii\,k\left(1+\log{(-\ii\pi)}\right)}{\pi}
=\int_{0}^{1}-\ii\,k(1-u)^{k-1}\cot{\frac{\pi u}{2}}+\frac{2\left(1-\pi^{-k}(\pi+\ii\log{(1-u)})^k\right)}{u^2}\,du
\end{multline}\\
\indent Now, another transformation is needed, for the purpose of separating the real and imaginary parts. Let us set $k=r+\ii\,t$, expressed in polar form, and change the integration variable using the relation $\log{(1-u)}=\pi\tan{v}$, chosen for convenience.
\begin{equation} \nonumber
k=r+\ii\,t=\sqrt{r^2+t^2}\exp\left(\ii\arctan{\frac{t}{r}}\right)
\end{equation}
\indent With that, taking into account the Jacobian of the transformation, the equation becomes,
\begin{multline} \nonumber
\frac{\pi(r-1)-2\,t(1+\log{\pi)}}{\pi}+\ii\frac{\pi\,t+2\,r(1+\log{\pi)}}{\pi}=\\ 
\pi\int_{0}^{\pi/2}(\sec{v})^2(-\ii\sqrt{r^2+t^2}\tan{\frac{\pi e^{-\pi\tan{v}}}{2}}\exp\left(-\pi\,r\tan{v}+\ii\left(\arctan{\frac{t}{r}}-\pi\,t\tan{v}\right)\right)\\+\frac{1}{2}\left(\csch{\frac{\pi\tan{v}}{2}}\right)^2\left(1-\exp\left(-r\log{\cos{v}}+t\,v+\ii(-t\log{\cos{v}}-r\,v)\right)\right))\,dv
\end{multline}\\
\indent Though this expression is very complicated, it can be simplified, as we do next. Since the parameters are real, the real and imaginary parts can be separated.

\subsection{The real part equation}
Below is the equation that is obtained for the real part,
\begin{multline} \nonumber
\frac{\pi(r-1)-2\,t(1+\log{\pi)}}{\pi}=\\ 
\pi\int_{0}^{\pi/2}(\sec{v})^2(\sqrt{r^2+t^2}\tan{\frac{\pi e^{-\pi\tan{v}}}{2}}\exp\left(-\pi\,r\tan{v}\right)\sin{\left(\arctan{\frac{t}{r}}-\pi\,t\tan{v}\right)}\\+\frac{1}{2}\left(\csch{\frac{\pi\tan{v}}{2}}\right)^2(1-\exp\left(-r\log{\cos{v}}+t\,v\right)\cos{(t\log{\cos{v}}+r\,v)}))\,dv
\end{multline}\\
\indent Any positive odd integer $r$ satisfies this equation, when $t=0$.\\

If $r+\ii\,t$ is a zero of the zeta function, so is its conjugate, $r-\ii\,t$. Hence, noting that the first term within the integral is an odd function in $t$, the above can be further simplified by adding the equations for $t$ and $-t$ as follows,
\begingroup
\small
\begin{multline} \nonumber
2(r-1)=\frac{\pi}{2}\int_{0}^{\pi/2}\left(\sec{v}\csch{\frac{\pi\tan{v}}{2}}\right)^2\left(2-(\sec{v})^r\left(e^{t\,v}\cos{(t\log{\cos{v}}+r\,v)}+e^{-t\,v}\cos{(t\log{\cos{v}}-r\,v)}\right)\right)\,dv
\end{multline}
\endgroup\\
\indent Let the function on the right-hand side of the equation be $f(r,t)$. After this transformation, the positive odd integers $r$ remain zeros of $f(r,0)=2(r-1)$.

\subsection{The imaginary part equation}
Likewise, for the imaginary part one has,
\begin{multline} \nonumber
\frac{\pi\,t+2\,r(1+\log{\pi)}}{\pi}=\\ 
\pi\int_{0}^{\pi/2}(\sec{v})^2(-\sqrt{r^2+t^2}\tan{\frac{\pi e^{-\pi\tan{v}}}{2}}\exp\left(-\pi\,r\tan{v}\right)\cos{\left(\arctan{\frac{t}{r}}-\pi\,t\tan{v}\right)}\\+\frac{1}{2}\left(\csch{\frac{\pi\tan{v}}{2}}\right)^2\exp\left(-r\log{\cos{v}}+t\,v\right)\sin{(t\log{\cos{v}}+r\,v)})\,dv
\end{multline}\\
\indent Coincidentally, any positive even integer $r$ satisfies this equation when $t=0$, so the two equations (real and imaginary) are never satisfied simultaneously for any positive integer.\\

Now the first term within the integral is an even function in $t$, so to simplify it the equations for $t$ and $-t$ can be subtracted, giving,
\begin{multline} \nonumber
2\,t=\frac{\pi}{2}\int_{0}^{\pi/2}(\sec{v})^{r+2}\left(\csch{\frac{\pi\tan{v}}{2}}\right)^2\left(e^{t\,v}\sin{(t\log{\cos{v}}+r\,v)}+e^{-t\,v}\sin{(t\log{\cos{v}}-r\,v)}\right)\,dv
\end{multline}\\
\indent Let the function on the right-hand side of the equation be $g(r,t)$. After this transformation, any $r$ satisfies $g(r,0)=0$, whereas the roots of $f(r,0)=2(r-1)$ are still the positive odd integers $r$. This means that when $t=0$, these transformations have introduced the positive odd integers $r$ as new zeros of the system, which were not there before.

\section{Riemann hypothesis reformulation}
If a linear combination of the equations for the real and imaginary parts is taken, such as $2(r-1)-2\,t\,\ii=f(r,t)-\ii\,g(r,t)$, the system of equations can be turned into a simpler single equation,
\begin{equation} \nonumber
k-1=\frac{\pi}{2}\int_{0}^{\pi/2}\left(\sec{v}\csch{\frac{\pi\tan{v}}{2}}\right)^2\left(1-\frac{\cos{k\,v}}{(\cos{v})^k}\right)\,dv
\end{equation}\\
\indent Going a little further, with a simple transformation ($u=\tan{v}$), the following theorem can be stated.\\

\textbf{Theorem} $k$ is a non-trivial zero of the Riemann zeta function if and only if $k$ is a non-trivial zero of,
\begin{equation} \label{eq:phi_integral}
\varphi(k)=1-k+\frac{\pi}{2}\int_{0}^{\infty}\left(\csch{\frac{\pi\,u}{2}}\right)^2\left(1-(1+u^2)^{k/2}\cos{(k\arctan{u})}\right)\,du
\end{equation}\\
\indent Hence the Riemann hypothesis is the statement that the zeros of  $\varphi(k)$ located on the critical strip have $\Re(k)=1/2$.\\

\textbf{Proof} All the roots of $\varphi(k)$ should also be roots of the zeta function, except for the positive odd integers and the trivial zeros of the eta function ($1+2\pi\ii\,j/\log{2}$, for any non-zero integer $j$), though this might not be true since the equations were transformed (that is, there might be other $k$ such that $\varphi(k)=0$ but $\zeta(k)\neq 0$).\\

From equations \eqrefe{eta_formula} and \eqrefe{h(k)_eta_rel}, the following relation between $h(k)$ and $\zeta(k)$ can be deduced,
\begin{equation}
h(k)=2\,\Gamma(k+1)(\ii\pi)^{-k}(1-2^{1-k})\zeta(k) \text{,}
\end{equation}
\noindent and from the observation that $\varphi(k)=\Re{(h(k))}$ for all real $k$, it follows that the relation between $\varphi(k)$ and $\zeta(k)$ is given by,
\begin{equation} \label{eq:FE}
\varphi(k)=-\frac{2\,\Gamma(k+1)\left(2^{1-k}-1\right)}{\pi^k}\cos{\frac{\pi k}{2}}\,\zeta(k) \text{,}
\end{equation}
\noindent for all complex $k$ except where undefined. This transformation takes the only pole of the zeta function to zero, hence $\varphi(k)$ is pole-free and entire (analytic on the whole complex plane).\footnote{Note that $\varphi(k)$ is not a Dirichlet series in any domain.} $\square$\\

Note this functional equation breaks down at the negative integers ($\Gamma(k+1)=\pm\infty$, but $\zeta(k)=0$ or $\cos{\pi k/2}=0$, whereas $0<\varphi(k)<\infty$) and at $k=1$. There is a trivial zeros trade-off between these two functions (the negative even integers for the positive odd integers). Note also that while the convergence domain of $h(k)$ is $\Re{(k)}>0$, the domain of $\varphi(k)$ is the whole complex plane.\\

If equation \eqrefe{FE} is combined with Riemann's functional equation, the following simpler equivalence is obtained, valid for all complex $k$ except the zeta pole,
\begin{equation} \nonumber
2(k-1)(1-2^{-k})\zeta(k)=\varphi(1-k) \text{}
\end{equation}

\subsection{Particular values of $\varphi(k)$ when $k$ is integer}
From the functional equations, the values of $\varphi(\pm k)$ when $k$ is a non-negative integer can be easily found out,
\begin{equation} \nonumber
\begin{cases}
  \varphi(2k)=\left(2-2^{2k}\right)B_{2k}\\
  \varphi(2k+1)=0\\ 
  \varphi(-k)=2k\left(1-2^{-k-1}\right)\zeta(k+1)
\end{cases}
\end{equation}\\
\indent The last identity implies that for large positive real $k$, $\varphi(-k)\sim 2k$. The function $\varphi(k)$ has trivial zeros at $k=1+2\pi\ii\,j/\log{2}$, where $j$ can be any integer, even zero. Those could be removed and the resulting function would still be entire.

\section{The $\varphi(k)$ generating function}
One can also create a generating function for $\varphi(k)$, based on the following identities,
\begin{equation} \nonumber
\cos{\arctan{u}}=\frac{1}{\sqrt{1+u^2}} \text{, and } \cos{(k\,\arctan{u})}=T_k(\cos{\arctan{u}}) \text{,}
\end{equation}
\noindent where $T_k(x)$ is the Chebyshev polynomial of the first kind.\\

Therefore, using the generating function of $T_k(x)$ available in the literature, for the positive integers one has,
\begin{equation} \nonumber
\sum_{k=0}^{\infty}\left(x\sqrt{1+u^2}\right)^k\,T_k(\cos{\arctan{u}})=\frac{1-x}{(1-x)^2+(x\,u)^2} \text{,}
\end{equation}
\noindent from which it is possible to produce the generating function of $\varphi(k)$,
\begin{equation} \label{eq:phi(x)_gf}
q(x)=\sum_{k=0}^{\infty}x^k\,\varphi(k)=\frac{2}{1-x}-\frac{1}{(1-x)^2}+\frac{\pi\,x^2}{2(1-x)}\int_{0}^{\infty}\left(\csch{\frac{\pi\,u}{2}}\right)^2\frac{u^2}{(1-x)^2+(x\,u)^2}\,du
\end{equation}\\
\indent The $k$-th derivative of $q(x)$ yields the value of $\varphi(k)$, and obtaining it is not very hard (we just  need to decompose the functions in $x$ into a sum of fractions whose denominators have degree one, if the roots are simple -- so we can easily generalize their $k$-th derivative). After we perform all the calculations and simplifications we find that,
\begin{equation} \nonumber
\varphi(k)=\frac{q^{(k)}(x)}{k!}=1-k-\frac{\pi\,k!}{2}\int_{0}^{\infty}\left(\csch{\frac{\pi\,u}{2}}\right)^2\sum_{j=1}^{\floor{k/2}}\frac{(-1)^j\,u^{2j}}{(2j)!(k-2j)!}\,du \text{,}
\end{equation}
\noindent where $\floor{k/2}$ means the integer division.\\

In here we went from an expression that holds for all $k$ to an expression that is only defined for $k$ a positive integer, loosely the opposite of analytic continuation.

\newpage

\section{Graphics plotting}
First we plot the curves obtained with the imaginary part equation, $g(r,t)$, with values of $r$ starting at $1/8$ with $1/8$ increments, up to $7/8$, for a total of 7 curves plus the $2\,t$ line. The points where the line crosses the curves are candidates for zeros of the Riemann zeta function (they also need to satisfy the real part equation, $f(r,t)=2(r-1)$).\\

Let us see what we obtain when we plot these curves with $t$ varying from 0 to 15. In the graph below, higher curves have greater $r$, though not always, below the $x$-axis it is vice-versa -- but generally the more outward the curve, the greater the $r$,
\begin{center}
\includegraphics{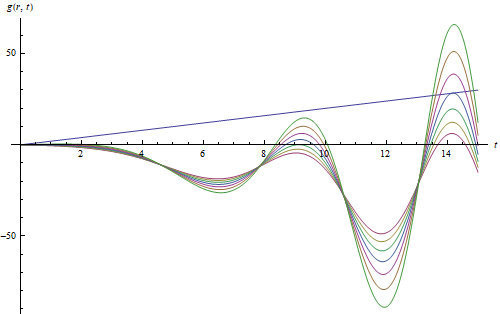}
\end{center}

As we can see, it seems the line crosses the curve for $r=1/2$ at its local maximum, which must be the first non-trivial zero (that is, its imaginary part). The line also crosses 3 other curves (all of which have $r>1/2$), but these are probably not zeros due to the real part equation. Also, it seems there must be a line that unites the local maximum points of all the curves, though that is just a wild guess.\\

One first conclusion is that one equation seems to be enough for $r=1/2$, the line seems to only cross this curve at the zeta zeros. Another conclusion is that apparently curves with $r<1/2$ do not even meet the first requirement, and also apparently $r=1/2$ is just right. A third conclusion is that all curves seem to have the same inflection points.\\

\newpage 

In the below graph we plotted $g(r,t)$ for the minimum, middle and maximum points of the critical strip (0, $1/2$ and 1), with $t$ varying from 0 to 26, for further comparison (0 is pink, 1 is green),
\begin{center}
\includegraphics{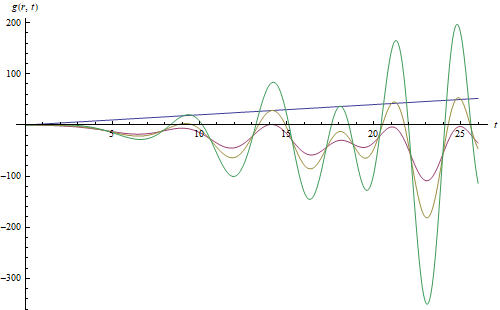}
\end{center}

\newpage 

Now, the graph below shows plots for curves $-2(r-1)+f(r,t)$ and $-2t+g(r,t)$ together. The plots were created for $r=0$ (red), $r=1/2$ (green) and $r=1$ (blue) (curves with the same color have the same $r$). A point is a zero of the zeta function when both curves cross the $x$-axis at the same point (three zeta non-trivial zeros are shown).
\begin{center}
\includegraphics{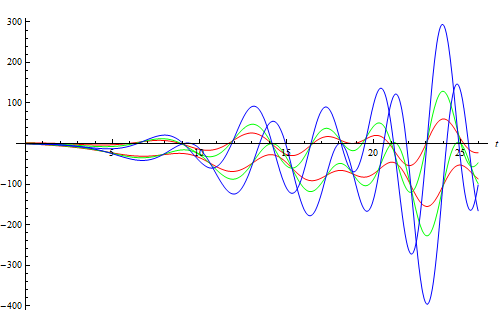}
\end{center}

\newpage 

And finally, graphs for the difference of the two functions, $-2(r-1)+f(r,t)+2t-g(r,t)$, were created for the same $r$'s as before and with the same colors as before (but now we also have pink ($r=1/4$) and cyan ($r=3/4$)). 
\begin{center}
\includegraphics{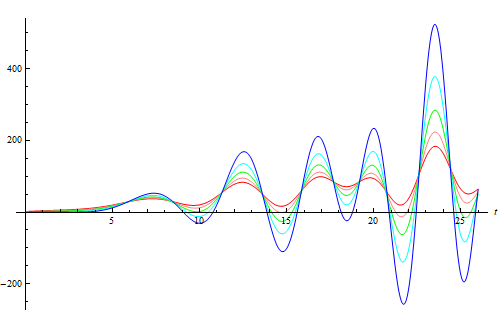}
\end{center}

\end{document}